\documentclass[12pt]{article}

\usepackage{lineno}
\nolinenumbers

\usepackage[utf8]{inputenc}
\usepackage[english]{babel}
\usepackage[table]{xcolor}
\usepackage[fleqn]{amsmath}
\usepackage{nccmath}
\usepackage{booktabs}
\usepackage{caption}
\usepackage{float}
\usepackage{sidecap}
\usepackage{csquotes}

\usepackage[backend=biber,style=alphabetic,sorting=ynt]{biblatex}

\usepackage{answers}
\usepackage{setspace}
\usepackage{graphicx}
\usepackage{enumitem}
\usepackage{multicol}
\usepackage{mathrsfs}
\usepackage[margin=1in]{geometry} 
\usepackage{amsmath,amsthm,amssymb}
\usepackage{mathtools}
\usepackage{adjustbox}

\newcommand{\Z}{\mathbb{Z}}
\newcommand{\C}{\mathbb{C}}

\DeclareRobustCommand{\Stirling}{\genfrac\{\}{0pt}{}}

\begin{document}

\title{A triangular field of rational numbers related to Stirling numbers and Hyperbolic functions}
\author{Andreas B. G. Blobel\\ %replace with your name
andreas.blobel@kabelmail.de} %if necessary, replace with your course title

\maketitle

\begin{abstract}
\noindent A triangular field of rational numbers is characterized, with relations to Stirling numbers $2^\text{nd}$ kind, Hyperbolic functions, and centered Binomial distribution. A Generating function is given.
\end{abstract}
\vspace{10mm}

\setcounter{secnumdepth}{0} % Suppress numbering of (sub)sections
\setcounter{tocdepth}{3} % Include 3 levels into TOC 

\tableofcontents

\subsection{Preliminary definitions}
%\addcontentsline{toc}{subsection}{Preliminary definitions}

Let $f(x)$ be a smooth function over some domain \emph{D}. Define the operator $[r]$ as follows:
\begin{ceqn}
\begin{align} \label{def_[r]}
    f^{[r]}(x) &:= \left[ x \cdot \tfrac{d}{dx} \right]^r f(x) & r  \in \mathbb{N}_0
\end{align}
\end{ceqn}

\vspace{5mm}\noindent
Please notice that this definition uses square brackets, differing from the familiar notation $f^{(r)}(x)$, which involves round brackets and indicates the $r^{\text{th}}$ derivative of $f(x)$.

\noindent
Define the elementary functions
%\begin{ceqn}
\begin{equation} \label{def_gh}
\begin{split}
    g,h  :\hspace{10pt} &\C \setminus \{0\}\ \longrightarrow\ \C \\[5pt]
    g(x) :=\; &\tfrac{1}{2}\,(x + x^{-1})      &h(x) := \tfrac{1}{2}\,(x - x^{-1})
\end{split}
\end{equation}
%\end{ceqn}
From \eqref{def_[r]} it follows immediately
\begin{equation} \label{gh[1]}
\begin{split}
    g^{[1]}(x) &= h(x)    & \hspace{49pt} h^{[1]}(x) = g(x)
\end{split}
\end{equation}
%
%\vspace{5mm}\noindent
Zeros of $g$ and $h$ are
\begin{equation} \label{Zeros_gh}
\begin{split}
    &g(i) \hspace{9pt} = 0   & \hspace{66pt} h(1) \hspace{9pt} = 0 \\
    &g(-i) = 0   &h(-1) = 0
\end{split}
\end{equation}
Here, $i$ denotes the imaginary unit.

\vspace{5mm}\noindent
Based on \eqref{def_gh}, define the functions
\begin{equation}  \label{def_Gsj}
\begin{split}
    G_{s,j} &:\hspace{9pt} \C \setminus \{0, 1, -1, i, -i\}\ \longrightarrow\ \C \\[8pt]
    G_{s,j} &:= g^{s-j} \cdot h^{j}  & \hspace{79pt} s \ge 0  \hspace{22pt} j \in \Z
\end{split}
\end{equation}
The domain of $G_{s,j}$, in general, excludes $\{0, 1, -1, i, -i\}$, because, if $j < 0$ the zeros of $h$
 \eqref{Zeros_gh} cause \emph{poles} \cite{Zeros_and_poles} of $G_{s,j}$, and so do the zeros of $g$ in case $j > s$.
For each given $s \ge 0$, the infinite set of functions
\begin{ceqn}
\begin{align} \label{independent}
    \mathcal{G}_s &= \big\{ \dots\ ,\ G_{s,-2}\ ,\ G_{s,-1}\ ,\ G_{s,\,0}\ ,\ G_{s,\,1}\ ,\ G_{s,\,2}\ ,\  \dots \big\} & s \ge 0
\end{align}
\end{ceqn}
are \emph{linearly independent}. The reason is, that for each $j \in \Z$, $G_{s,j}$ has \emph{unique orders} of zeros/poles. For example, by definition \eqref{def_Gsj}, $x = 1$ is a zero of order $j$ if $j > 0$, and a pole of order $j$ if $j < 0$. Therefore, no finite linear combination of \emph{other} elements in $\mathcal{G}_s$ can reproduce the same order $j$ of the zero/pole at $x = 1$.

\subsection{The triangular array $A$}

Application of the  $[r]$ operator with $r=1$ to $G_{s,j}$ as defined in \eqref{def_Gsj} yields
\begin{align} \label{G[1]}
    G^{[1]}_{s,j} &= j \cdot G_{s,j-1} + (s-j) \cdot G_{s,j+1} & s \ge 0  \hspace{22pt} j \in \Z
\end{align}
This suggests the expansion
\begin{align} \label{G[r]_expan}
     G^{[r]}_{s,0} &= \sum_{j \in \Z} A_{s,r,j} \cdot G_{s,j} & s,r \ge 0
\end{align}
Applying $[r = 1]$ to \eqref{G[r]_expan} leads to a recurrence relation for the coefficients:

\vspace{5mm} \noindent
Here are the steps in detail:
\begin{fleqn} % this is for left align; \usepackage{nccmath} needed in preamble 
\begin{equation}\label{ArecurrDeriv}
\begin{alignedat}{5}
&G^{[r+1]}_{s,0} &&= \sum_{j \in \mathbb{Z}} A_{s,r,j} &&\cdot G^{[1]}_{s,j} \\[5pt]
&                &&= \sum_{j \in \mathbb{Z}} A_{s,r,j} \cdot j &&\cdot G_{s,j-1}
                 &&\;+\; \sum_{j \in \mathbb{Z}} A_{s,r,j} \cdot (s-j) &&\cdot G_{s,j+1} \\[5pt]
&                &&= \sum_{j \in \mathbb{Z}} A_{s,r,j+1} \cdot (j+1) &&\cdot G_{s,j}
                 &&\;+\; \sum_{j \in \mathbb{Z}} A_{s,r,j-1} \cdot (s-(j-1)) &&\cdot G_{s,j} \\[5pt]
\end{alignedat}
\end{equation}
\end{fleqn}
\eqref{A_recurr} follows, if one compares this with \eqref{G[r]_expan} after having replaced $r \to r+1$, and when observing linear independence  of \eqref{independent}.
%\vspace{8pt}\noindent
\begin{subequations}
\begin{equation} \label{A_recurr}
\begin{split} 
     A_{s,r+1,j} \;&=\; A_{s,r,j-1} \cdot (s - (j-1)) \\[5pt]
                   &\,+\; A_{s,r,j+1} \cdot (j+1) & \hspace{62pt} j \in \Z \hspace{22pt} s,r \ge 0
\end{split}
\vspace{4pt}
\end{equation}
\begin{align}
    A_{s,0,\,0} \hspace{9pt} &=\; 1  \label{A_init} \\
    A_{s,r,j} \hspace{11pt} &=\; 0 \hspace{27pt} \textsc{if not} \hspace{5pt} 0 \le j \le r  \label{A_triang}
\end{align}
\end{subequations}
Initial or \emph{seed values} \eqref{A_init} / \eqref{A_triang}
make $A_{s,r,j}$ a \emph{Triangular array} \cite{Triangular_array}. We can therefore write \eqref{G[r]_expan} as a \emph{finite} sum:
\begin{align} \label{G[r]_expanFinite}
     G^{[r]}_{s,0} &= \sum_{j =0}^{r} A_{s,r,j} \cdot G_{s,j} & s,r \ge 0
\end{align}

\vspace{5pt}\noindent
The coefficients in \eqref{G[1]} add up to $j + (s-j) = s$.
By induction, this is readily extended to the identity
\begin{align} \label{A_sumFinite}
    \sum_{j =0}^{r} A_{s,r,j} &= s^r & s,r \ge 0
\end{align}

\vspace{5pt}\noindent
Table \ref{table:A} displays a few elements from the tip of $A_{s,r,j}$, as computed from \eqref{A_recurr}. Here, the $(s)_k$ denote \emph{Falling factorials} \cite{Falling_factorials}. 
The elements of $A$ are integer functions of the size parameter $s \ge 0$.
In particular, elements in the $j=0$ column can be expressed in terms of $\cosh$ \eqref{A0cosh}. At the same time, for any given $r \ge 0$, the $r^\text{th}$ item in the $j=0$ column equals the $r^\text{th}$ moment of a centered Binomial distribution of size $s$ \eqref{mu_r_M0}.
Coefficients of $(s)_k$ in the $j=0$ column match the triangle of numbers given in \cite{A156289}, when selecting only even values of $r$.

%\vspace{5mm} \noindent
\begin{table}[ht!]
\begin{adjustbox}{width=1\textwidth}
\centering
\begin{tabular}{r r | l l l l l} 
  & & $j$\\
  & & 0 & 1 & 2 & 3 & 4 \\
 \hline
 $r$ & 0 & $1$ & $\cdot$ & $\cdot$ & $\cdot$ & $\cdot$ \\[5pt]
 & 1 & 0 & $s$     & $\cdot$ & $\cdot$ & $\cdot$ \\[5pt]
 & 2 & $s$ & 0 & $(s)_2$ & $\cdot$ & $\cdot$ \\[8pt]
 & 3 & 0 & $3\,(s)_2 + s$ & 0 & $(s)_3$ & $\cdot$ \\[8pt]
 & 4 & $3\,(s)_2 + s$ & 0 & $6\,(s)_3 + 4\,(s)_2$ & 0 & $(s)_4$ \\[8pt]
 & 5 & 0 & $15\,(s)_3 + 15\,(s)_2 + s$ & 0 & $10\,(s)_4 + 10\,(s)_3$ & 0 \\[8pt]
 & 6 & {$15\,(s)_3 + 15\,(s)_2 + s$} & 0 & $45\,(s)_4 + 75\,(s)_3 + 16\,(s)_2$
 & $0$ & $\dots$ \\[8pt]
 & 7 & 0 & $105\,(s)_4 + 210\,(s)_3 + 63\,(s)_2 + s$ & $\dots$ \\[8pt]
 & 8 & $105\,(s)_4 + 210\,(s)_3 + 63\,(s)_2 + s$ & $\dots$ \\[8pt]
 & 9 & 0 & $\ddots$ \\[8pt]

  %\hline
\end{tabular}
\end{adjustbox}
\caption{Apex of $A_{s,r,j}$}
\label{table:A}
\end{table}

\subsubsection{Relation to Hyperbolic functions}

Exactly the same table is generated if one replaces $g\,/\,h$ from \eqref{def_gh} by $\cosh/\sinh$ \cite{Hyperbolic_functions}, and, at the same time, replaces the $[r]$ operator from \eqref{def_[r]} by the $r^\text{th}$ derivative $(r)$:
\begin{ceqn}
\begin{equation} \label{cosh_corr}
    \begin{split}
        g\; &\longrightarrow\; \cosh \\
        h\; &\longrightarrow\; \sinh \\
        [r]\; &\longrightarrow\; (r) \\
    \end{split}
\end{equation}
\end{ceqn}
With these replacements \eqref{G[r]_expanFinite} becomes
\begin{align} \label{cosh^s(r)_expanFinite}
     \Big(\cosh(t)^s\Big)^{(r)}\; &= \sum_{j =0}^{r} A_{s,r,j}\, \cdot\, \cosh(t)^{s-j} \, \sinh(t)^j & s,r \ge 0
\end{align}
and therefore, setting the formal parameter $t = 0$:
\begin{align} \label{A0cosh}
    A_{s,r,\,0} &\;=\; \Bigg[\; \Big(\cosh(t)^s\Big)^{(r)}\; \Bigg]_{t=0}  & s,r \ge 0
\end{align}

\vspace{5mm} \noindent
\subsubsection{Relation to centered Binomial distribution}

$G_{s,0}$ as defined in \eqref{def_Gsj}, if written as the \emph{Laurent} series
\begin{align} \label{Gs_Laurent}
     G_{s,0}(x) &= \sum_{j \in \mathbb{Z}} b_{s,j} \cdot x^j & s \ge 0
\end{align}
generates a centered \emph{Binomial} distribution \cite{Binomial_distribution}. For if we define the set of \emph{events}
\begin{align} \label{Vset}
    \mathcal{V}_s &:= \big\{ 2k-s\; |\; 0 \le k \le s \big\} & s\ge 0
\end{align}
each event $v \in \mathcal{V}_s$ is associated with $b_{s,v}$ from \eqref{Gs_Laurent}, which can be written as
\begin{align} \label{CBD} %Centered Binomial Distribution
    b_{s,v} &= \frac{1}{2^s}\, \binom{s}{\tfrac{v+s}{2}} & v \in \mathcal{V}_s\;,\; s \ge 0
\end{align}
According to e.g. \cite{Bulmer} the \emph{Moment generating function} $M_s(t)$ is defined by the \emph{Expected value} of $e^{t\, \cdot\, v}$
\begin{align} \label{MGF}
    M_s(t) &= \sum_{k=0}^s \frac{1}{2^s} \binom{s}{k} \; e^{t\, \cdot\, (2k-s)} \;=\; \cosh(t)^s   & s \ge 0
\end{align}
Therefore, if $\mu_{s,r}$ denotes the $r^\text{th}$ moment of \eqref{CBD} and $M_s^{(r)}$ denotes the $r^\text{th}$ derivative of $M_s(t)$ with respect to the formal parameter $t$ , we have
\begin{align} \label{mu_r_M0}
    \mu_{s,r} &= M_s^{(r)}(0) \;=\; \Bigg[\; \Big(\cosh(t)^s\Big)^{(r)}\; \Bigg]_{t=0}  & s,r \ge 0
\end{align}

\vspace{5mm} \noindent
On the other hand, by construction of $G^{[r]}_{s,0}$, we have the identity
\begin{align} \label{mu_r_G1}
    \mu_{s,r} &= G^{[r]}_{s,0}(1) \;=\; \sum_{j = 0}^{r}\; A_{s,r,j}\, \cdot\, G_{s,j}(1) \;=\; A_{s,r,\,0}  & s,r \ge 0
\end{align}
Comparing \eqref{mu_r_M0} and \eqref{mu_r_G1} again gives \eqref{A0cosh}.

\vspace{8pt}
\subsubsection{The condensed form $B$}
In order to exclude the intermediate zero elements in table \ref{table:A}, we derive from $A$ the ``condensed'' form
\begin{align} \label{B_def}
     B_{s,r,n} \hspace{7pt} :&= \hspace{2pt} A_{s,r,r-2n}    & \hspace{34pt} n \in \Z \hspace{22pt}  s,r \ge 0
\end{align}
Transcribing \eqref{A_recurr}, \eqref{A_init} / \eqref{A_triang}, and \eqref{A_sumFinite} in terms of \eqref{B_def} gives:
\begin{subequations}
\begin{equation} \label{B_recurr}
\begin{split} 
     B_{s,r+1,n} \;=\; &B_{s,r,\,n} \cdot (s - (r - 2n)) \\[5pt]
                +\; &B_{s,r,\,n-1} \cdot (r - 2(n-1)) & \hspace{108pt} s,r \ge 0
\end{split}
\end{equation}
%\vspace{3pt}
\begin{align}
 B_{s,0,\,0} \hspace{9pt} &=\; 1   \label{B_init}\\
 B_{s,r,n} \hspace{11pt} &=\; 0 \hspace{27pt}  \textsc{if not} \hspace{5pt} 0 \le n \le \big\lfloor \tfrac{r}{2} \big\rfloor  \label{B_triang}
\end{align}
\end{subequations}

\vspace{9pt}
\begin{align} \label{B_sumFinite}
    \sum_{n = 0}^{\big\lfloor \tfrac{r}{2} \big\rfloor}\; B_{s,r,n} &= s^r  & s,r \ge 0
\end{align}

\vspace{1cm}
\subsection{The triangular field $\varphi$}

\vspace{5mm}
\subsubsection{Expansion in terms of falling factorials}
Consider the transformation
\begin{align} \label{phi_def}
    B_{s,r,n} &= \sum_{0 \le k \le j \le n} (s)_{r - n - j} \cdot (r)_{2n+k} \cdot \varphi_{n,j,k}
    & s,r \ge 0 \hspace{22pt} n \in \Z 
\end{align}

\vspace{3mm}\noindent
If one inserts \eqref{phi_def} into \eqref{B_recurr}, while applying appropriate identities to match factorial powers \cite{Falling_factorials} of $s$ and $r$ on both sides,
one gets the recurrence relation
\begin{fleqn}
\begin{subequations}
\begin{equation}\label{phi_recurr}
\begin{split}
(2n + k) \cdot \varphi_{n,j,\,k}\; =\; (n - (j-1))\; \cdot\; & \varphi_{n,j-1,\,k-1} \\[5pt]
+\; (k+1)\;\; \cdot\; & \varphi_{n-1,j,\,k+1} \\[5pt]
+\; &\varphi_{n-1,j,\,k} & \hspace{57pt} 0 \le k \le j \le n
\end{split}
\end{equation}
\vspace{5pt}
\begin{align}
\varphi_{0,0,0}\; &= \; 1  \label{phi_init} \\[5pt]
\varphi_{n,j,k}\; &= \; 0  & \hspace{93pt} \textsc{if not} \hspace{8pt} 0 \le k \le j \le n  \label{phi_triang}
\end{align}
\end{subequations}
\end{fleqn}

\vspace{11pt} \noindent
\emph{Seed} values \eqref{phi_init} / \eqref{phi_triang}, which are counterparts of conditions \eqref{B_init} / \eqref{B_triang}, make $\varphi$ a triangular array in 3 dimensions. For given $n \ge 0$, $\varphi_{n,j,k}$ can be regarded as a \emph{lower triangular matrix} of dimension $n+1$. Indices $j$ and $k$ denote row and column indices respectively. The matrix associated with $n = 0$ reduces to the scalar $1$.

\vspace{6pt}
\subsubsection{Special solutions}
Some special solutions of the triple  \eqref{phi_init} / \eqref{phi_triang} / \eqref{phi_recurr} are easily verified:
\begin{subequations}
\begin{align}
    \varphi_{n,k,k} \hspace{10pt} &= \frac{1}{2^n\,(n-k)!\;k!\;3^k} & n \ge k \ge 0 \label{phi_nkk}\\[7pt]
    \varphi_{n,n-1,0} &= \frac{1}{(2n)!} \hspace{49pt} \varphi_{n,n,1}= \frac{1}{(2n+1)!} & n \ge 1 \label{phi_lambda=1}\\[9pt]
    \varphi_{n,n,0} \hspace{10pt} &= \delta_{n,0} & n \in \Z \label{phi_nn0}
\end{align}
\end{subequations}

\vspace{5mm} \noindent
From \eqref{phi_nkk} we have in particular
\begin{align} \label{phi_n00nnn}
    \varphi_{n,0,0} \hspace{11pt} &= \frac{1}{2^n\,n!} \hspace{49pt} \varphi_{n,n,n}= \frac{1}{6^n\,n!} & \hspace{44pt} n \ge 0
\end{align}

\vspace{5mm} \noindent
\subsubsection{Relation to Stirling numbers}
Inserting \eqref{phi_def} into \eqref{B_sumFinite} gives, after rearrangement
\begin{align} \label{B_sum_phi}
    \sum_{i=0}^{r}\; (s)_{r-i} \; \sum_{j \ge 0} \; (r)_j \; \Bigg[ \sum_{n \ge 0}\; \varphi_{n,i-n,j-2n} \Bigg]   &= s^r  & s,r \ge 0
\end{align}
 Comparing this with the basic relation for Stirling numbers $2^\text{nd}$ kind \cite{Stirling_second}
\begin{align} \label{Stirling_sum}
    \sum_{i = 0}^{r} (s)_{r-i} \; \Stirling{r}{r-i} &= s^r & s,r \ge 0
\end{align}
we get the identity
\begin{align} \label{phi_Stirling}
     \sum_{j \ge 0}\; (r)_j \; \Bigg[ \sum_{n \ge 0}\; \varphi_{n,i-n,j-2n} \Bigg]  &= \Stirling{r}{r-i}  & 0 \le i \le r
\end{align}

\vspace{5mm} \noindent
\subsubsection{Relation to Hyperbolic functions}
If one substitutes $r = 2n$ in \eqref{phi_def} one gets
\begin{equation}
\begin{split} \label{B2n}
    B_{s,2n,n} &= \sum_{0 \le k \le j \le n} (s)_{n - j} \cdot (2n)_{2n+k} \cdot \varphi_{n,j,k} \\[8pt]
               &=\; (2n)_{2n} \cdot \sum_{0 \le j \le n} (s)_{n - j} \cdot \varphi_{n,j,\,0}
    & \hspace{91pt} s,n \ge 0
\end{split}
\end{equation}
And therefore, from \eqref{B_def} and \eqref{A0cosh}:
\begin{align} \label{phi2cosh}
    \sum_{0 \le j \le n} (s)_{n - j} \cdot \varphi_{n,j,\,0} \;&= \; \frac{1}{(2n)!}\; \Bigg[\; \Big(\cosh(t)^s\Big)^{(2n)}\; \Bigg]_{t=0}
    &  s,n \ge 0
\end{align}

\vspace{6pt}
\subsubsection{The adjoint form $\Tilde{\varphi}$}
It is useful to introduce the \emph{adjoint} form
\begin{align} \label{phiT_def}
\Tilde{\varphi}_{n, \lambda,\, k} &:=   \varphi_{n, n - \lambda + k,\, k} & \lambda,k \in \Z \hspace{38pt} n \ge 0
\end{align}
Rewriting  \eqref{phi_recurr} / \eqref{phi_init} / \eqref{phi_triang} in terms of $\Tilde{\varphi}$ gives
\begin{fleqn}
\begin{subequations}
\begin{equation}\label{phiT_recurr}
\begin{split}
(2n + k) \cdot \Tilde{\varphi}_{n,\lambda,\,k}\; =\; (\lambda - (k-1))\; \cdot\;\;
& \Tilde{\varphi}_{n,\lambda,\,k-1} \\[5pt]
+\; (k+1)\;\; \cdot\;\; & \Tilde{\varphi}_{n-1,\lambda,\,k+1} \\[5pt]
+\; & \Tilde{\varphi}_{n-1,\lambda-1,\,k} & \hspace{59pt}  0 \le k \le \lambda \le n
\end{split}
\end{equation}
\vspace{5pt}
\begin{align}
\Tilde{\varphi}_{0,0,0}\; &= \; 1 \label{phiT_init} \\[5pt]
\Tilde{\varphi}_{n,\lambda,\,k}\; &= \; 0  & \hspace{105pt} \textsc{if not} \hspace{8pt} 0 \le k \le \lambda \le n  \label{phiT_triang}
\end{align}
\end{subequations}
\end{fleqn}

\noindent
The \textasciitilde modifier in \eqref{phiT_def} can be viewed as a mapping which turns the lower triangular matrix $\varphi_n$ into the lower triangular matrix $\Tilde{\varphi}_{n}$. In that sense, it is an \emph{involution}:
\begin{align} \label{involution}
\Tilde{\Tilde{\varphi}}_{n} &=   \varphi_{n} & n \ge 0
\end{align}

%\pagebreak
\begin{table}[ht!]
%\begin{center}
\begin{adjustbox}{width=1.0\textwidth}
%\begin{tabular}{ l | l l l l l l l l l l} 
\begin{tabular}{ l | c c c c c c c c c c} 
 & $-1$ & $0$ & $1$ & $2$ & $\cdots$ & $\lambda-1$ & $\lambda$ & $\cdots$ & $n$ & $n+1$ \\
\hline
$-1$ & \multicolumn{10}{c}{\cellcolor{cyan}}  \\[3mm]
$0$ & \cellcolor{cyan} & $\frac{1}{2^n\,n!}$ & \multicolumn{8}{c}{\cellcolor{cyan}}  \\[3mm]
$\vdots$ & \cellcolor{cyan}  & $\vdots$ & $\ddots$ & \multicolumn{7}{c}{\cellcolor{cyan}}  \\[3mm]
$n-\lambda$ & \cellcolor{cyan}  & $\Tilde{\varphi}_{n,\lambda,\,0}$ & & & \multicolumn{6}{c}{\cellcolor{cyan}} \\[3mm]
$n-\lambda+1$ & \cellcolor{cyan} & & $\Tilde{\varphi}_{n,\lambda,\,1}$
& & & \multicolumn{5}{c}{\cellcolor{cyan}}    \\[3mm]
$n-\lambda+2$ & \cellcolor{cyan}  & & & $\Tilde{\varphi}_{n,\lambda,\,2}$
& & & \multicolumn{4}{c}{\cellcolor{cyan}} \\[3mm]
$\vdots$  & \cellcolor{cyan}  &  & & & $\ddots$ & & & \multicolumn{3}{c}{\cellcolor{cyan}} \\[3mm]
$n-1$  & \cellcolor{cyan} & $\frac{1}{(2n)!}$
& & & $\cdots$ & $\Tilde{\varphi}_{n,\lambda,\,\lambda-1}$ & & $\ddots$ & \multicolumn{2}{c}{\cellcolor{cyan}} \\[3mm]
$n$  & \cellcolor{cyan} & $\delta_{n,0}$  & $\frac{1}{(2n+1)!}$
& & $\cdots$ & & $\Tilde{\varphi}_{n,\lambda,\,\lambda}$ & $\cdots$ & $\frac{1}{6^n\,n!}$ & \cellcolor{cyan} \\[3mm]
$n+1$  & \multicolumn{10}{c}{\cellcolor{cyan}}  \\[3mm]
\end{tabular}
\end{adjustbox}
%\end{center}
\caption{Triangular field $\varphi_{n,j,k}$ / $\Tilde{\varphi}_{n,\lambda,\,k}$}
\label{table:phi}
\end{table}

\vspace{5mm}\noindent
Table \ref{table:phi} illustrates the triangular shape of $\varphi_{n,j,k}$ / $\Tilde{\varphi}_{n,\lambda,\,k}$. Row and column indices, $j$ and $k$, run over all integers $\Z$. The blue area, which extends in all directions, marks zero elements. The \emph{diagonals} of  $\varphi_{n,j,k}$ form the \emph{rows} of the adjoint field $\Tilde{\varphi}_{n,\lambda,\,k}$ \eqref{phiT_def}, and vice versa, by involution \eqref{involution},
the \emph{diagonals} of $\Tilde{\varphi}_{n,\lambda,\,k}$ form the \emph{rows} of $\varphi_{n,j,k}$.
Special solutions \eqref{phi_n00nnn}, \eqref{phi_lambda=1}, and \eqref{phi_nn0} have been inserted.

\vspace{1cm}\noindent
\subsection{The Generating functions $F_{n,\lambda}$}
Based upon the adjoint quantities \eqref{phiT_def}, define the functions
\begin{align} \label{Fnl_def}
    F_{n,\lambda}(z) :&= \sum_{k \in \Z}\ \Tilde{\varphi}_{n,\lambda,\,k} \cdot z^k 
    & \hspace{21pt} n, \lambda \in \Z
\end{align}
Here, $z$ is a formal parameter. If $F^{'}_{n,\lambda}$ denotes the derivative of \eqref{Fnl_def} with respect to $z$, we have
\begin{align} \label{Fnl_deriv}
    F^{'}_{n,\lambda}(z) &= \sum_{k \in \Z}\ k \cdot \varphi_{n,n-\lambda+k,\, k} \cdot z^{k-1}  & n, \lambda \in \Z
\end{align}
If one multiplies \eqref{phiT_recurr} by $z^k$, performs summation over all $k \in \Z$, matches powers of $z$ on both sides, and replaces \eqref{Fnl_def} and \eqref{Fnl_deriv}, one gets the \emph{recursive differential equation}
\begin{fleqn}
\begin{subequations}
\begin{equation}  \label{Fnl_differential}
    \begin{split}
(2n - \lambda\,z)\ &F_{n,\lambda}(z)\; +\; z(z+1)\ F^{'}_{n,\lambda}(z) \hspace{14pt} = \\[5pt]
&F_{n-1,\lambda - 1}(z)\hspace{25pt} +\ F^{'}_{n-1,\lambda}(z)
& \hspace{44pt} n \ge \lambda \ge 0
    \end{split}
\end{equation}
\vspace{0pt}
\begin{equation} \label{Fnl_init}
\hspace{55pt} F_{0,\,0} \; = \; 1
\end{equation}
%\vspace{0pt}
\begin{equation} \label{Fnl_triang}
\hspace{55pt} F_{n,\lambda} \; = \; 0  \hspace{123pt} \textsc{if not} \hspace{8pt} n \ge \lambda \ge 0
\end{equation}
\end{subequations}
\end{fleqn}
with boundary conditions \eqref{Fnl_init} and \eqref{Fnl_triang} being direct consequences of the corresponding properties \eqref{phiT_init} and \eqref{phiT_triang} of $\Tilde{\varphi}$. A distinctive special value is
\begin{align} \label{Fnl0}
    F_{n,\lambda}(0) &= \Tilde{\varphi}_{n,\lambda,\,0} =  \varphi_{n,n-\lambda,\,0}
    & \hspace{21pt} n, \lambda \in \Z
\end{align}

\vspace{5mm}
\subsubsection{Frequency representation of integer partitions}
The solution \eqref{Fnl_solution} of \eqref{Fnl_differential} involves integer partitions in  \emph{frequency representation} \cite{Frequency_Representation}.
Let $\mathcal{P}_{n,\lambda}$ denote the set of integer partitions of $n \ge 0$ with a given number $\lambda \ge 0$ of parts \cite{IntegerPartition}.
Let $\pi \in \mathcal{P}_{n,\lambda}$ be a partition of $n$ with $\lambda$ parts. Then $\pi(k)$ denotes the \emph{frequency} of part $k \ge 1$ in $\pi$. This implies the identities
\begin{ceqn}
\begin{subequations}
\begin{align}
    n &= \sum_{k = 1}^{n}\ \pi(k) \cdot k \label{pi_sum_parts} \\
    \lambda &= \sum_{k = 1}^{n}\ \pi(k) \label{pi_sum_freq}
\end{align}
\end{subequations}
\end{ceqn}
In the degenerate case $\lambda = 0$ we have
\begin{ceqn}
\begin{equation} \label{EmptyPartition}
\mathcal{P}_{n,\ 0} = \begin{cases}
\big\{ \textit{empty partition} \big\} &\text{if\; $n = 0$}\\[8pt]
\hspace{40pt} \emptyset &\text{if\; $n > 0$}
\end{cases}
\end{equation}
\end{ceqn}
\eqref{EmptyPartition} means, that $n = 0$ has exactly one partition, namely the \emph{empty partition}, without parts. And on the other hand, clearly, there is no partition of $n > 0$ without parts.

\vspace{5mm}
\subsubsection{Solution of \eqref{Fnl_differential} / \eqref{Fnl_init}}
 The recursive differential equation \eqref{Fnl_differential}, including boundary condition \eqref{Fnl_init}, is solved by
\begin{align} \label{Fnl_solution}
F_{n,\lambda}(z) \;&=\; \sum_{\pi \in \mathcal{P}_{n,\lambda}}\ \prod_{k=1}^{n}
\Bigg[ \frac{1}{\pi(k)!} \cdot f_k(z)^{\pi(k)} \Bigg]
& n \ge \lambda \ge 0
\end{align}
the component functions $f_k$ being defined as
\begin{align} \label{Fk_def}
f_{k}(z) \;:&=\;  { \frac{1}{(2k)!}\ \bigg(1\, +\, \frac{z}{2k+1} } \bigg)
& k \ge 0 
\end{align}

\vspace{3mm}\noindent
Setting $z=0$ in \eqref{Fnl_solution} while using \eqref{Fnl0} gives
\begin{align} \label{Fnl0_solution}
\varphi_{n,n-\lambda,\,0} \;&=\; \sum_{\pi \in \mathcal{P}_{n,\lambda}}\ 
\Bigg[ \prod_{k=1}^{n}\ \pi(k)! \cdot \big( (2k)! \big)^{\pi(k)} \Bigg]^{-1}
& n \ge \lambda \ge 0
\end{align}
For $\lambda=0$ this reduces to \eqref{phi_nn0}, and for $\lambda=n$ to the first part of \eqref{phi_n00nnn}.

\vspace{1cm}
\subsection{Proof that \eqref{Fnl_solution} solves \eqref{Fnl_differential} / \eqref{Fnl_init}}
The proof is divided into three sections.
First, get rid of the recursiveness in \eqref{Fnl_differential} by forming power series over the whole range of the indices $\lambda$ and $n$. This leads to the linear partial differential equation \eqref{F_differential} including boundary condition \eqref{F_init}.
Second, solve \eqref{F_differential} / \eqref{F_init}.
And third, focus on specific powers of $x$ and $y$ in the representation of solution \eqref{F_solution} as an infinite product of infinite sums \eqref{F_factors}.

\vspace{7pt}
\subsubsection{Summation over $\lambda \in \Z$}
Starting from \eqref{Fnl_def}, define the functions
\begin{align} \label{Fn_def}
    F_{n}(y,z) :&= \sum_{\lambda \in \Z}\ F_{n,\lambda}(z) \cdot y^\lambda 
    & \hspace{21pt} n \in \Z
\end{align}
Here, $y$ and $z$ are formal parameters.
If one multiplies \eqref{Fnl_differential} by $y^\lambda$, performs summation over all $\lambda \in \Z$, matches powers of $y$ and $z$ on both sides, and replaces \eqref{Fn_def}, or their partial derivatives, one gets the \emph{recursive} partial differential equation
\begin{fleqn}
\begin{subequations}
\begin{equation}  \label{Fn_differential}
\begin{split}
2\,n\, &F_{n}(y,z) \;+\; z(z+1)\ \partial_z F_{n}(y,z) \;-\; z\,y\, \partial_y F_{n}(y,z)\hspace{8pt} =\\
   y\, &F_{n-1}(y,z) \hspace{37pt} +\; \partial_z F_{n-1}(y,z) \hspace{133pt}  n \ge 0
\end{split}    
\end{equation}
\vspace{0pt}
\begin{equation} \label{Fn_init}
\hspace{17pt} F_{0} \; = \; 1
\end{equation}
%\vspace{0pt}
\begin{equation} \label{Fn_triang}
\hspace{17pt} F_{n} \; = \; 0  \hspace{260pt} n < 0
\end{equation}
\end{subequations}
\end{fleqn}
with boundary conditions \eqref{Fn_init} and \eqref{Fn_triang} again being direct consequences of the corresponding properties \eqref{Fnl_init} and \eqref{Fnl_triang}. A distinctive special value here is
\begin{align} \label{Fn00}
    F_{n}(0,0) \ &=\  F_{n,\,0}(0)\ =\ \Tilde{\varphi}_{n,\,0,\,0}\ =\  \varphi_{n,n,\,0}\ =\ \delta_{n,0}
    & \hspace{11pt} n \in \Z
\end{align}
Here, we have made use of \eqref{Fnl0} and \eqref{phi_nn0}.

\vspace{7pt}
\subsubsection{Summation over $n \in \Z$}
Next, in a quite analogous manner, define the function
\begin{align} \label{F_def}
    F(x,y,z) :&= \sum_{n \in \Z}\ F_{n}(y,z) \cdot x^n
\end{align}
where $F_{n}$ has been declared in \eqref{Fn_def}.
Again, $x$, $y$, and $z$ are formal parameters.
If one multiplies \eqref{Fn_differential} by $x^n$, performs summation over all $n \in \Z$, matches powers of $x$, $y$, and $z$ on both sides, and replaces \eqref{F_def}, or their partial derivatives, one gets the linear partial differential equation
\begin{subequations}
\begin{multline} \label{F_differential}
2\,x\, \partial_x F(x,y,z) \,+\, z(z+1)\, \partial_z F(x,y,z) \,-\, z\,y\, \partial_y F(x,y,z)\ =\\
y\,x\, F(x,y,z) \,+\, x\, \partial_z F(x,y,z)
\end{multline}
%\vspace{0pt}
\begin{equation} \label{F_init}
\hspace{27pt} F(0,0,0) \ = \ F_{0}(0,0) \ =\ 1
\end{equation}
\end{subequations}
with boundary condition \eqref{F_init} derived from \eqref{Fn00}.

\vspace{7pt}
\subsubsection{Solution of \eqref{F_differential}}
The linear partial differential equation \eqref{F_differential}, including boundary condition \eqref{F_init}, is solved by
 \begin{ceqn}
\begin{equation} \label{F_solution}
F(x,y,z) \;=\; \exp{
\Bigg[\ y \cdot \sum_{k \ge 1}\, f_k(z) \cdot x^k \ \Bigg] }
\end{equation}
 \end{ceqn}
 where $f_k$ is from \eqref{Fk_def}. This is verified by insertion. Result \eqref{Fnl_solution} becomes clear
if one writes down \eqref{F_solution} as an infinite product of infinite sums, as illustrated in \eqref{F_factors}. Here, each row represents an infinite sum.

\vspace{7pt}
\begin{fleqn} % this is for left align; \usepackage{nccmath} needed in preamble 
\begin{equation}\label{F_factors}
\begin{alignedat}{5}
&1 \hspace{18pt}+\hspace{18pt}                 y   \ f_1(z)   \ &&x^{1\,\cdot\,1}
   \hspace{18pt}+\hspace{18pt} \tfrac{1}{2!} \ y^2 \ f_1^2(z) \ &&x^{1\,\cdot\,2}
   \hspace{18pt}+\hspace{18pt} \tfrac{1}{3!} \ y^3 \ f_1^3(z) \ &&x^{1\,\cdot\,3}
   \hspace{18pt}+\hspace{18pt} &&\cdots\\[12pt]
&1 \hspace{18pt}+\hspace{18pt}                 y   \ f_2(z)   \ &&x^{2\,\cdot\,1}
   \hspace{18pt}+\hspace{18pt} \tfrac{1}{2!} \ y^2 \ f_2^2(z) \ &&x^{2\,\cdot\,2}
   \hspace{18pt}+\hspace{18pt} \tfrac{1}{3!} \ y^3 \ f_2^3(z) \ &&x^{2\,\cdot\,3}
   \hspace{18pt}+\hspace{18pt} &&\cdots\\[12pt]
&1 \hspace{18pt}+\hspace{18pt}                 y   \ f_3(z)   \ &&x^{3\,\cdot\,1}
   \hspace{18pt}+\hspace{18pt} \tfrac{1}{2!} \ y^2 \ f_3^2(z) \ &&x^{3\,\cdot\,2}
   \hspace{18pt}+\hspace{18pt} \tfrac{1}{3!} \ y^3 \ f_3^3(z) \ &&x^{3\,\cdot\,3}
   \hspace{18pt}+\hspace{18pt} &&\cdots\\[12pt]
&\vdots &&\vdots &&\vdots &&\vdots &&\ddots
\end{alignedat}
\end{equation}
\end{fleqn}
\eqref{Fnl_solution} follows for any given pair $n, \lambda \ge 0$, if, upon multiplication of all rows, one combines all coefficients of $x^n$, and at the same time, all coefficients of $y^\lambda$.

\vspace{1cm}
\subsection{Summary}
Starting from the pair of elementary functions \eqref{def_gh}, the triangular array $A$ (see Table \ref{table:A}) has been generated, whose elements are integer functions of the size parameter $s \ge 0$.
The 3-dimensional triangular field of rational numbers $\varphi$ has been introduced, being a result of the transformation \eqref{phi_def}, which involves the condensed form $B$ \eqref{B_def}.
$\varphi$ satisfies the triple of conditions \eqref{phi_recurr} / \eqref{phi_init} / \eqref{phi_triang}.
 Distinctive special solutions have been given in \eqref{phi_nkk}, \eqref{phi_lambda=1}, \eqref{phi_nn0}, and \eqref{Fnl0_solution}. Relations to Stirling numbers \eqref{phi_Stirling} and Hyperbolic functions \eqref{phi2cosh} have been described. A Generating function \eqref{Fnl_def}, based upon the adjoint form $\Tilde{\varphi}$ \eqref{phiT_def}, has been defined, and it was shown that it can be written as \eqref{Fnl_solution}.
%\medskip

%\pagebreak
\listoftables
\printbibliography

\end{document}